\documentclass[12pt]{article}
\usepackage{amssymb,amsmath}
\begin{document}
\newtheorem{thm}{Theorem}[section]
\newtheorem{lem}[thm]{Lemma}

\title{\textbf{The Skitovich-Darmois theorem for discrete and compact totally disconnected Abelian groups}}
\author{I.P. Mazur}
\maketitle

\section{Introduction}

The classical Skitovich-Darmois theorem states (\cite{Skit},\cite{Darm}, see also \cite[Ch. 3]{KLR}):
Let $\xi_i ,i=1,2,\ldots,n,$ $n\geq2,$ be independent random variables, and $\alpha_i,\beta_i$ be nonzero numbers. Suppose that the linear forms $L_1=\alpha_1\xi_1+\cdots+\alpha_n\xi_n$ и $L_2=\beta_1\xi_1+\cdots+\beta_n\xi_n$ are independent. Then all random variables $\xi_i $ are Gaussian.

The Skitovich-Darmois theorem was generalized into various algebraic structures, in particular, into locally compact Abelian groups(\cite{F book}-\cite{FG 2005},\cite{Gr F},\cite{Mazur}).  In these researches the random variables take values in a locally compact Abelian group $X$, coefficients of the linear forms are topological automorphisms of $X$, and the number of the linear forms are two.
In the article we continue these researches and study the Skitovich-Darmois theorem in the case when random variables take values in different classes of locally compact Abelian groups but the number of linear forms more than 2.

Throughout the article $X$ is a second countable locally compact Abelian group.
Let $Aut(X)$ be the group of topological automorphisms of $X$, $\mathbb{Z}(k)=\{0,1,2\ldots,k-1\}$
be the group of residue modulo  $k$. Let $x\in X$.
Denote by $E_x$ the degenerate distribution, concentrated at the point $x$. If $K$ is a compact subgroup of  $X$,
denote by $m_K$ the Haar distribution on $K$.
Denote by $I(X)$ the set of shifts of such distributions, i.e. the distributions of the form $m_K\ast E_x$,
where $K$ is a compact subgroup of $X$, $x\in X$.

The first step to generalize the Skitovich-Darmois theorem on locally compact Abelian groups is to consider the case of finite Abelian groups. This problem was studied first in \cite{F 1992}. Note that the idempotent distributions on a finite Abelian group can be regarded as analogues of the Gaussian distributions on real line. Let $\xi_i, i=1,2,\ldots,n,n\geq 2,$ be independent random variables with values in a finite Abelian group $X$
and distributions $\mu_i$. Let $\alpha_j,\beta_j$ be automorphisms of $X$. Consider the linear forms $L_1=\alpha_1\xi_1+\cdots+\alpha_n\xi_n$ и $L_2=\beta_1\xi_1+\cdots+\beta_n\xi_n$.
G.M. Feldman proved that
the class of groups, on which the independence of $L_1$ и $L_2$ implies that all $\mu_i\in I(X)$, are poor and consists of the groups of the form  (\cite{F 1992})
\begin{equation}
\mathbb{Z}(2^{m_1})\times\cdots\times\mathbb{Z}(2^{m_l}),0\leq m_1<\cdots<m_l. \label{Lestnica}
\end{equation}
On the other hand if we consider $n$ linear forms of $n$ independent random variables, then the Skitovich-Darmois theorem is valid for an arbitrary finite Ableian group. Namely, the following theorem holds (\cite{Mazur}):

\begin{thm}\label{Th 1}
Let $\xi_i, i=1,2,\ldots,n,n\geq 2,$ be independent random variables with values in a finite Abelian group $X$ and distributions $\mu_i$.
If the linear forms $L_j=\sum_{i=1}^{n}\alpha_{ij}\xi_i,$
where $\alpha_{ij}\in Aut(X),i,j=1,2,\ldots,n$, are independent, then $\mu_i\in I(X),i=1,2,\ldots,n$.
\end{thm}

Our first aim is to generalize Theorem \ref{Th 1} into discrete groups. Using this result we extend Theorem \ref{Th 1} into groups of the form $\mathbb{R}^m\times D \times K$, where $D$ is a discrete, $K$ is a compact totally disconnected group of the special form, $m\geq 0$. Finally we describe conditions on compact totally disconnected group, which are necessary and sufficient for the validity of the Theorem \ref{Th 1}.

\section{Definitions and designations}

Let $X$ be a second countable locally compact Abelian group.
Let $\{X_{\lambda}:\lambda\in \Lambda\}$ be a non empty family of groups. Denote by $\mathbf{P}_{\lambda\in \Lambda}X_{\lambda}$ the direct product of the groups $X_{\lambda}$. If $X=X_{\lambda}$ for all $\lambda\in \Lambda $, then the direct product of the groups $X_{\lambda}$ we denote by $X^{\mathfrak{n}}$, where $\mathfrak{n}$ is a cardinal number of $\Lambda$.

Let $p$ be a prime number. Denote by $\mathfrak{P}$ the set of all prime number.
Denote by $\triangle_p$ the group of $p$-adic integers. We recall that an Abelian group is called a $p$-primary if every nonzero element of this group has order $p^m$ for some natural $m$.

Let $\mu$ be a distribution on $X$. The characteristic function of the distribution $\mu$ is defined by the formula (\cite{F book APD})
$$\hat{\mu}(y)=\int_X(x,y)d\mu(x),y\in Y.$$

Recall (\cite{Part}), that a distribution $\mu$ on $X$ is called Gaussian, if its characteristic function can be represented in the form
$$\hat{\mu}(y)=(x,y)\exp\{-\phi(y)\}, \quad y\in Y,$$
where the function $\phi(y)$ satisfy the equation $$\phi(u+v)+\phi(u-v)=2\phi(u)+2\phi(v),\quad u,v\in Y.$$
Denote by $\Gamma(X)$ the set of all Gaussian distributions on $X$.

\section{Main results}

The main result of the article is Theorem \ref{Theorem Prod Gr}.
The prof of Theorem \ref{Theorem Prod Gr} uses the following theorem, which is interesting in itself.

\begin{thm}
\label{Theorem DAG} Let $X$ be a discrete Abelian group. Let $\xi_i, i=1,2,\ldots,n,$ be independent random variables with values in $X$ and distributions $\mu_i$.
Consider the linear forms $L_j=\sum_{i=1}^{n}\alpha_{ij}\xi_i,$
where $\alpha_{ij}\in Aut(X),i,j=1,2,\ldots,n$.
Then the independence of the linear forms $L_j$ implies that $\mu_i=E_{x_i}\ast m_{G_i}$, where $G_i$ are finite subgroups of $X$, $x_i\in X$, $i=1,2,\ldots,n$.
\end{thm}

Note that Theorem \ref{Theorem DAG} was proved earlier in (\cite[13.17]{F book}) in the case $n=2$.

\begin{thm}
\label{Theorem Prod Gr} Let $X$ be a Abelian group of the form
\begin{equation}\label{th prod X=}
X=\mathbb{R}^m\times K\times D,
\end{equation}
where, $m\geq 0$ $K$ --- compact totally disconnected group, which is topologically isomorphic to a group of the form
\begin{equation}
\quad \mathbf{P}_{p \in \mathfrak{P}}(\triangle_p^{n_p}\times G_p),\label{Komp Vp nesv gr}
\end{equation}
where $n_p$ are some non negative integers, $G_p$ are finite $p$-primary groups, а $D$ is a discrete countable group.
Let $\xi_i, i=1,2,\ldots,n,$ be independent random variables with values in $X$ and distributions $\mu_i$.
Then the independence of the linear forms $L_j=\sum_{i=1}^{n}\alpha_{ij}\xi_i,$
where $\alpha_{ij}\in Aut(X),i,j=1,2,\ldots,n,$ implies that $\mu_i\in \Gamma(X)\ast I(X), i=1,2,\ldots,n$.
\end{thm}

From the conditions of theorem {\ref{Theorem Prod Gr}} it follows that condition (\ref{Komp Vp nesv gr})
 is sufficient for validity of the Skitovich-Darmois theorem on compact totally disconnected Abelian groups
 in the case of $n$ linear forms of $n$ random variables. It turns out, that this condition is necessary.

\begin{thm}
Let $X$ be a  compact totally disconnected Abelian group, $\xi_i,i=1,2,\ldots,n,$be independent random variables with values in $X$ and distributions $\mu_i$, $\alpha_{ij}\in Aut(X)$. From the independence of the linear forms $L_j=\sum_{i=1}^{n}\alpha_{ij}\xi_i,j=1,2,\ldots,n,$ follow that all $\mu_i \in \Gamma(X)\ast I(X)$, if and only if  $X$ is topologically isomorphic to the group of the form $(\ref{Komp Vp nesv gr})$.
\end{thm}

Note that the theorem above was proved earlier in (\cite[13.25]{F book}) in the case $n=2$.

\textbf{Remark.} There are locally compact Abelian groups which are non isomorphic to the group of the form (\ref{th prod X=}), and on which The Skitovich-Darmois theorem is valid in the case of $n$ linear forms of $n$ random variables.
The needed example was constructed for the case $n=2$ in \cite[раздел 3]{FG 2010}. It easily can be checked, that this example can be generalized in the case of the arbitrary $n$.

Mathematical Division\\
B. Verkin Institute for Low \\
Temperature Physics and Engineering\\
 of the National Academy \\
 of Sciences of Ukraine\\
47, Lenin Ave, Kharkov\\
61103, Ukraine

\end{document}